\documentclass[12pt]{article}
\usepackage[utf8]{inputenc}
\usepackage[russian]{babel}
\usepackage{amssymb}
\usepackage{amsmath}

\begin{document}

\begin{center}
\textbf{\large Change of variable and the rapidity\\ of decrease of Fourier
coefficients}
\end{center}

\begin{center}
\textsc{Vladimir Lebedev}
\end{center}

\begin{quotation}
{\small \textbf{Abstract.} We consider the class $C(T)$ of
continuous real-valued functions on the circle. For certain
classes of functions naturally cha\-ra\-cte\-rised by the
rapidity of decrease of Fourier coefficients we inve\-sti\-gate
whether it is possible to bring families of functions in $C(T)$
into these classes by a change of variable. This paper was
originally published in \emph{Matematicheski\v{\i} Sbornik},
\textbf{181}:8 (1990), 1099--1113 (Russian). The English
trans\-la\-ti\-on, published in \emph{Mathematics of the USSR,
Sbornik}, \textbf{70}:2 (1991), 541--555, is to a large extent
inconsistent with the original text. Herein the author provides
a corrected trans\-la\-tion.


  }
\end{quotation}

\begin{center}
\textbf{Introduction}
\end{center}

  We consider the class $C( T)$ of continuous real-valued functions
on the circle $ T$, and the Fourier series of functions in $C(
T)$:
$$
f\sim\sum_{k\in  Z}\widehat{f}(k)e^{ikt}
\eqno(1)
$$
($ Z$ is the set of integers).

  It is known that certain properties of functions of class
$C(T)$, related to the series (1), can be improved by a change
of variable, i.e., by a homeomorphism of $T$ onto itself. For a
survey of basic results in the area and a number of open
problems see the papers by A. M. Olevski\v{\i} [1], [2].
Following [1] and [2] we quote some of these results.

\quad

\textsc{Theorem} A (J. P\'al, 1914; H. Bohr, 1935). \emph{For
every function $f\in C( T)$ there exists a homeomorphism $h$ of
$ T$ onto itself such that the superposition $f\circ h$ belongs
to the class $U(T)$ of functions with uniformly convergent
Fourier series.}

\quad

  The method used to prove this theorem allows us to obtain the rapid decrease
of the Fourier coefficients, namely
$$
|\widehat{f\circ h}(k)|=d(k)+O(1/|k|), \qquad |k|\rightarrow\infty,
$$
where
$$
\sum_{k\in Z}|d(k)|^2|k|<\infty,
$$
whence, in particular, we see that
$$
\sum_k |\widehat{f\circ h}(k)|^p<\infty \qquad \forall p>1.
\eqno(2)
$$

  On the other hand, we have the following theorem that solves a problem
posed by N. N. Lusin.

\quad

\textsc{Theorem} B (A. M. Olevski\v{\i}, 1981). \emph{There
exists an $f\in C(T)$ such that there is no change of variable
which will bring $f$ into the algebra $A(T)$ of absolutely
convergent Fourier series, i.e., $f\circ h\notin A(T)$ for
every homeomorphism $h$ of $T$ onto itself.}

\quad

Thus, in general, it is impossible to attain condition (2) for
$p=1$.

On the other hand, the following theorem shows that one can
approach the class $A(T)$ arbitrarily close.

\quad

\textsc{Theorem} C (A. A. Saakyan, 1979). \emph{ If $\alpha(n),
n=0, 1, 2, \ldots,$ is a positive sequence satisfying the
condition $\sum_n \alpha(n)=\infty$ and a certain condition of
regularity, then for every $f\in C(\mathbb T)$ there is a
homeomorphism $h$ such that $\widehat{f\circ
h}(k)=O(\alpha(|k|))$} (see [1], Theorem 4.1).

\quad

J.-P. Kahane and Y. Katznelson investigated whether it is
possible to bring families of functions into the class $U(T)$.
They obtained the following result (1978).

\quad

\textsc{Theorem} D. \emph{For every compact set $K$ in $C(T)$
there is a change of variable which brings $K$ into $U(T)$,
i.e., there exists a homeomorphism $h : T\rightarrow T$ such
that $f\circ h\in U( T)$ for all $f\in K$.}

\quad

  The same authors, considered the classes
$$
A_\varepsilon(T)=\bigg\{x : \|x\|_{A_\varepsilon}
=\sum_{k\in Z}|\widehat{x}(k)|\varepsilon_{|k|}<\infty\bigg\},
$$
where $\varepsilon=\{\varepsilon_n\}$ is a sequence of positive
numbers that tends to zero, and proved (in 1981) the following
theorem:

\quad

\textsc{Theorem} E. \emph{There exists a sequence
$\varepsilon_n\rightarrow 0$ and a pair of functions in $C(T)$
such that there is no change of variables which will bring the
pair into $A_\varepsilon(T)$.}

\quad

  We notice that, if $\varepsilon_n\rightarrow 0$, then for every
individually taken function there is a change of variable which
brings it into $A_\varepsilon$ (this follows from Theorem C).

  In this paper we consider the classes
$A_\varepsilon(T)$ and certain other classes of functions
naturally characterised by the rate of decrease of Fourier
coefficients and investigate if it is possible to bring compact
families of functions in $C(T)$ into these classes.

  We show (\S~2) that \emph{under certain assumptions of
regularity of the sequence $\varepsilon$ the condition
$$
\sum_{n=1}^\infty \varepsilon_n/n<\infty
$$
is necessary and sufficient in order that for every compact set
$K\subset C(T)$ there is a change of variable which brings $K$
into $A_\varepsilon$.}

Here the sufficiency is a direct consequence of the following
assertion (\S~1): \emph{if $\lambda_n\rightarrow\infty$ and $K$
is a compact subset of $C(T)$, then there exists a
homeomorphism $h :
 T \rightarrow T$ such that
$$
|\widehat{f\circ h}(k)|=O(\lambda_{|k|}/|k|)  \qquad \forall f\in K.
$$
}

It is not clear if one can put $\lambda_n\equiv 1$ in this
assertion. However, it is impossible to attain the condition
$|\widehat{f\circ h}(k)|=o(1/|k|), ~\forall f\in K$. We show
(\S~1) that \emph{there exists a compact set such that there is
no change of variable which will bring it into the class $\{x :
|\widehat{x}(k)|=o(1/|k|)\}$.} This gives an answer to the
problem posed in [1] (Russian p. 182, English p. 210) and in
[2], \S~4.2.

  It follows from Theorem C that, for every function $f\in C(T)$,
some superposition $x=f\circ h$ with a homeomorphism $h : T
\rightarrow T$ satisfies
$$
\sum |\widehat{x}(k)|^2|k|<\infty.
\eqno(3)
$$
We show (\S~3) that even for two functions it is, in general,
impossible to attain (3) by a single homeomorphism: \emph{if
$f\in C(T)$ has unbounded variation on $T$, then there exists
$g\in C(T)$ such that there is no change of variable which will
bring the pair $\{f, g\}$ into the class defined by condition
\emph{(3)}.} Thus the answer to the problem posed in [3], p.
41, is negative. On the other hand, \emph{if $f\in C( T)$ is a
function of bounded variation, then every pair $\{f, g\}, ~g\in
C( T),$ can be brought to the indicated class.}

  We shall use notation $H^\omega( T)$ for the class of functions $f$ on
$T$ that satisfy $\omega(\delta, f)=O(\omega(\delta))$, where
$$
\omega(\delta, f)=\sup_{|t_1-t_2|<\delta}|f(t_1)-f(t_2)|
$$
is the uniform modulus of continuity of $f$, and
$\omega(\delta)$ is a given increasing continuous function on
$[0, \infty)$ with $\omega(0)=0$.

\begin{center}
\textbf{\S~1. Estimates for $|\widehat{f\circ h}(k)|, ~f\in
H^\omega$}
\end{center}

\textsc{Theorem 1.} \emph{Let $\lambda_n\rightarrow\infty$ as
$n\rightarrow\infty$. Then, given any $\omega$, there exists a
homeomorphism $h$ of $ T$ onto itself such that
$$
|\widehat{f\circ h}(k)|=O(\lambda_{|k|}/|k|), \qquad |k|\rightarrow\infty
$$
for every $f\in H^\omega( T)$.}

\quad

  We define the integral modulus of continuity of a summable function
$x$ (of class $L( T)$) on $ T$ by
$$
\omega_1(\delta, x)=
\sup_{|\varepsilon|<\delta}\|x(\cdot+\varepsilon)-x(\cdot)\|_{L( T)}.
$$

  Theorem 1 is an immediate consequence of the well-known estimate
$|\widehat{x}(k)|=O(\omega_1(1/|k|, x))$ and the following
assertion:

\quad

\textsc{Lemma 1.} \emph{Let $\lambda(\delta)\rightarrow 0$ as
$\delta\rightarrow 0$. Then, given any $\omega$, there exists a
homeomorphism $h$ of $ T$ onto itself such that
$$
\omega_1(\delta, f\circ h)= O(\lambda(\delta)\delta),
\qquad \delta\rightarrow 0,
\eqno(4)
$$
for every $f\in H^\omega(T)$.}

\quad

\textsc{Proof.} We use a modification of the method used to
prove Theorem D (see [1], Theorem 4.2).

  For an arbitrary set $E\subset T$ we denote its $\delta$
-neighbourhood by $(E)_\delta$. Let $L_\infty(T)$ denote the
space of essentially bounded functions on $T$. Let $|E|$ denote
the Lebesgue measure of a set $E$.

 Let us prove first a simple lemma.

 \quad

 \textsc{Lemma 2.} \emph{Let $E\subset
T$ be a closed set. Suppose that $x\in L_\infty(T)$ is a
function which is constant on each interval complementary to
$E$. Then
$$
\omega_1(\delta, x)\leq 2\|x\|_{L_\infty}|(E)_\delta|.
$$
}

\textsc{Proof.} For $\varepsilon>0$ we have
$$
\int_ T |x(t+\varepsilon)-x(t)| dt=
\int_{(E)_\varepsilon} |x(t+\varepsilon)-x(t)| dt
\leq 2 \|x\|_{L_\infty}|(E)_\varepsilon|.
$$
Therefore $\omega_1(\delta, x)\leq 2\|x\|_{L_\infty}|(E)_\delta|$.

\quad

   We now proceed to the proof of Lemma 1. Without
loss of generality we may assume that $\lambda(\delta)$ is
monotonic and that $\liminf\lambda(\delta)\delta=0$. Consider a
nowhere dense perfect set $E\subset[0, 2\pi]$ that contains the
points $0$ and $2\pi$ and satisfies $|(E)_\delta|\leq
\lambda(\delta)\delta, ~\forall \delta>0$. For an interval
$I=(a, b)\subset [0, 2\pi]$ by $E(I)$ we denote the image of
$E$ under a homothetic mapping of $[0, 2\pi]$ onto $[a, b]$. It
is easy to verify that
$$
|(E(I))_\delta|\leq \lambda(\delta)\delta \qquad \forall \delta>0.
\eqno(5)
$$

    Let $N_k, k=1, 2, \ldots,$ be a sequence of positive integers with
$N_1=1$. Each such sequence defines sets $E_{k l}, ~k=1, 2,
\ldots, ~l=1, \ldots, N_k,$ as follows. We put $E_{1 1}=E=E((0,
2\pi))$. If the sets $E_{k l}, ~k=1, 2, \ldots, r, ~l=1,
\ldots, N_k$ have already been defined, we put $E_{r+1
l}=E(I_{r+1 l}), ~l=1, \ldots, N_{r+1,}$ where $I_{r+1 l},
~l=1, 2, \ldots,$ are the intervals complementary to
$$
\bigcup_{k=1}^r\bigcup_{l=1}^{N_k} E_{k l},
$$
enumerated in the order of nonincreasing length.

  It is clear that if the numbers $N_k$ increase fast enough, then
$$
\sup_l |I_{r l}|\rightarrow 0, \qquad r\rightarrow\infty.
\eqno(6)
$$

  We choose a sequence $\{N_k\}$ and define sets $E_{k l}, ~k=1, 2,
\ldots, ~l=1, \ldots, N_k,$ so that (6) is satisfied.

  Now we notice that if $g$ is a continuous function on $[0, 2\pi]$, then
we have an expansion
$$
g\overset{L_\infty}{=}\sum_{k=1}^\infty\sum_{l=1}^{N_k}g_{k l},
\eqno(7)
$$
where each function $g_{k l}$ is constant on the intervals
complementary to the corresponding $E_{k l}$, the series
converges in the $L_\infty$ norm, and, in addition,
$$
\|g_{1 1}\|_{L_\infty}\leq\|g\|_{L_\infty}, \qquad
\|g_{k l}\|_{L_\infty}\leq |g(I_{k l})|, \quad l=1, \ldots, N_k, ~k>1.
$$
To prove this we proceed as follows. Consider a set of the form
$E(I)$ and a function $x$ continuous on $I$. Let $P_{E(I)}(x)$
stands for a function that takes constant value $(x(a)+x(b))/2$
on each interval $(a, b)\subset I$ complementary to $E(I)$ and
vanishes at the other points of $[0, 2\pi]$. We put $g_{1
1}=P_{E_{1 1}}(g)$. If the functions $g_{k l}, ~l=1, \ldots,
N_k, ~k=1, \ldots, r$ have already been constructed, we put
$$
g_{r+1 l}=P_{E_{r+1 l}}\bigg(g-\sum_{k=1}^r\sum_{l=1}^{N_k}g_{k l}\bigg),
\qquad l=1, \ldots, N_{r+1}.
$$
Continuing this process, we obtain (7). Indeed,
$$
\bigg|\bigcup_{k=1}^\infty\bigcup_{l=1}^{N_k}E_{k l}\bigg|=0
$$
and, taking (6) into account, we have
$$
\bigg\|g-\sum_{k=1}^r\sum_{l=1}^{N_k}g_{k l}\bigg\|_{L_\infty}\leq
\sup_l \omega(|I_{r+1 l}|, g)\rightarrow 0.
$$
The rest of the properties of $g_{k l}$'s are obvious.

  Now, for each set $E_{k l}$ we fix a continuous increasing
function $h_{k l}$ that is constant on the intervals
complementary to $E_{k l}$ in $(0, 2\pi)$ and satisfies $h_{k
l}(0)=0$ and $h_{k l}(2\pi)=2\pi$. Let
$$
h=\sum_{k=1}^\infty\delta_k\sum_{l=1}^{N_k}h_{k l}.
$$
We choose the numbers $\delta_k ~k=1, 2 \ldots,$ so that
$$
\sum_{k=1}^\infty \delta_k N_k=1,
\eqno(8)
$$
and
$$
\sum_{k=1}^\infty N_k\omega(\varepsilon_k)<\infty,
\eqno(9)
$$
where
$$
\varepsilon_k=2\pi\sum_{s=k}^\infty \delta_s N_s.
$$
It follows from (6) and (8) that $h$ is a homeomorphism $
T\rightarrow T$.

   Note that
$$
|h(I_{k l})|\leq
\sum_{s=k}^\infty \delta_s\sum_{l=1}^{N_s}|h_{s l}(I_{k l})|\leq
2\pi \sum_{s=k}^\infty \delta_s N_s=\varepsilon_k,
\quad l=1, \ldots, N_k, \quad k=1, 2, \ldots.
$$

  Let $f\in H^\omega( T)$. Then for the function $g=f\circ h$
we have expansion (7), where
$$
\|g_{1 1}\|_{L_\infty}\leq \|g\|_{L_\infty},\qquad \|g_{k l}\|_{L_\infty}
\leq |f\circ h(I_{k l})|\leq \omega(\varepsilon_k, f),
\quad l=1, \ldots, N_k, ~k>1.
$$
Each function $g_{k l}$ is constant on the intervals
complementary to the corresponding set $E_{k l}$, so, by Lemma
2, taking (5) into account, we obtain
$$
\omega_1(\delta, g_{1 1})\leq 2\|g\|_{L_\infty}\lambda(\delta)\delta,
$$
$$
\omega_1(\delta, g_{k l})
\leq 2\omega(\varepsilon_k, f)\lambda(\delta)\delta \quad \forall\delta>0,
~k>1, ~l=1,\ldots, N_k,
$$
which, together with (9), proves Lemma 1. Theorem 1 follows.

\quad

  We note that for an individually taken function in $C(
T)$ estimate (4)  was established by B. S. Kashin (see [1],
Russian p. 179, English p. 206).

   Theorem 1 implies the following corollary:

\quad

\textsc{Corollary 1.} \emph{For every class $H^\omega( T)$
there exists a change of variable which brings it into
$\bigcap_{p>1}A_p( T)$, where $A_p( T)$ is the class of
functions on $T$ defined by
$$
x\in A_p( T)\Leftrightarrow
\sum_{k\in Z}|\widehat{x}(k)|^p<\infty.
$$}

\quad

  This result was stated by A. M. Olevski\v{\i} in [1] (Russian
p. 182, English, p. 210).

  The following question related to Theorem 1 is open: Can we
attain the estimate
$$
|\widehat{f\circ h}(k)|=O(1/|k|) \quad \forall f\in H^\omega(T)?
$$
It is also not clear whether it is possible to attain such an
estimate for an arbitrary pair of functions in $C(T)$.

  From Theorem C it is clear that if $f\in C(T)$, then
for a certain ho\-me\-o\-mor\-phism $h :  T\rightarrow T$ we
have
$$
|\widehat{f\circ h}(k)|=o(1/|k|).
\eqno(10)
$$
Is it true that for every class $H^\omega(T)$ there is a
ho\-me\-o\-mor\-phism $h$ such that (10) holds for all $f\in
H^\omega(T)$? This question was posed in [1] and [2]. Theorem 2
below gives a negative answer to this question. A similar
question for pairs of functions is open.

\quad

\textsc{Theorem 2.} \emph{There exists a class $H^\omega( T)$
such that for every ho\-me\-o\-mor\-phism $h$ of $ T$ onto
itself, the condition
$$
|\widehat{f\circ h}(k)|=o(1/|k|) \quad \forall f\in H^\omega( T)
\eqno(11)
$$
fails to be satisfied.}

\quad

\textsc{Proof.} Let $F\subset T$ be a perfect set. Suppose that
$E$ is a set of uniqueness. It is known (see [4], Chapter XIV,
\S~11) that if a function $g$ is constant on the intervals
complementary to $F$ and satisfies
$$
|\widehat{g}(k)|=o(1/|k|), \quad |k|\rightarrow\infty,
$$
then $g$ is equivalent to (i.e., coincides almost everywhere
with) a constant function.

  Let $\omega(\delta)=(\log (1/\delta))^{-\alpha}$,
where $0<\alpha<1$. Let $h$ be a homeomorphism $ T\rightarrow
T$. We shall construct a function $f\in H^{\omega}( T),
~f\not\equiv\mathrm{const},$ such that the superposition
$f\circ h$ is constant on the intervals complementary to a
certain perfect set which is a set of uniqueness, thus the
theorem will follow. Recall a result of Kahane and Salem (see
[5], Chapter VII, \S~8): if
$$
\liminf_{\varepsilon\rightarrow 0}
\frac{N_\varepsilon(F)}{\log(1/\varepsilon)}=0,
\eqno(12)
$$
where $N_\varepsilon(F)$ is the smallest number of intervals of
length $\varepsilon$ that cover $F$, then $F$ is a set of
uniqueness.

  Let $n_s$ be the positive integer closest to
$2^{2^{s/\alpha}}, ~s=1, 2, \ldots$. Obviously we have
$$
2^k/\log \prod_{s=1}^k n_s\rightarrow 0, \qquad k\rightarrow\infty,
\eqno(13)
$$
$$
2^k\omega\bigg(2\pi\prod_{s=1}^k\frac{1}{2n_s+1}\bigg)\geq \gamma>0,
\qquad k=1, 2, \ldots,
\eqno(14)
$$
where $\gamma$ is independent of $k$.

  Denote the homeomorphism inverse to $h$ by $h^{-1}$.
We partition the interval $[-\pi, \pi]$ into $2n_1+1$ pairwise
nonoverlapping closed intervals of equal length. Enumerating
them in order of succession on $[-\pi, \pi]$, we find, among
the first $n_1$ intervals, the one, which we denote by $I_0$,
such that
$$
|h^{-1}(I_0)|\leq 2\pi/n_1.
$$
Similarly, among the last $n_1$ intervals we find an interval
$I_1$ such that
$$
|h^{-1}(I_1)|\leq 2\pi/n_1.
$$
Suppose the closed intervals $I_{i_1\ldots i_s}, ~s=1, \ldots,
k,$ where $i_s=0$ or $1$, have already been constructed. We
partition $I_{i_1\ldots i_k}$ into $2n_{k+1}+1$ closed
intervals of equal length. Enumerating them in order of
succession on $I_{i_1\ldots i_k}$, we find, among the first
$n_{k+1}$ and the last $n_{k+1}$ intervals, the intervals which
we denote by $I_{i_1\ldots i_k 0}$ and $I_{i_1\ldots i_k 1}$
respectively, such that
$$
|h^{-1}(I_{i_1\ldots i_ki_{k+1}})|\leq
\frac{1}{n_{k+1}}|h^{-1}(I_{i_1\ldots i_k})|, \qquad i_{k+1}\in\{0; 1\}.
$$

  Thus we have a correspondence between tuples
$(i_1, \ldots, i_k), ~k=1, 2, \ldots,$ of zeros and ones, and
closed intervals $I_{i_1\ldots i_k}, ~k=1, 2, \ldots,$ with the
following properties:
$$
I_{i_1\ldots i_ki_{k+1}}\subset I_{i_1\ldots i_k},
\qquad |h^{-1}(I_{i_1\ldots i_k})|\leq 2\pi\bigg/\prod_{s=1}^k n_s,
\eqno(15)
$$
$$
\inf\{|t_0-t_1|: t_0\in I_{i_1\ldots i_k 0}, ~t_1\in I_{i_1\ldots i_k 1}\}
\geq \delta_{k+1}=
$$
$$
=2\pi\bigg/\prod_{s=1}^{k+1}(2n_s+1), \qquad k=1, 2, \ldots.
\eqno(16)
$$

Let
$$
E_k=\bigcup_{\alpha\in \{0; 1\}^k} I_\alpha
$$
(the union is taken over all $k$ -tuples of zeros and ones). We
put $E=\bigcap_{k=1}^\infty E_k$. Obviously, $E$ is a nowhere
dense perfect set. It follows from (15) that the set
$F=h^{-1}(E)$ can be covered by $2^k$ intervals of length
$2\pi/\prod_{s=1}^k n_s, ~k=1, 2, \ldots$. Taking (13) into
account, we see that $F$ satisfies (12), so $h^{-1}(E)$ is a
set of uniqueness. To complete the proof of the theorem it
remains to construct a function $f\in H^\omega( T),
~f\not\equiv\mathrm{const},$ which is constant on the intervals
complementary to $E$.

  For each closed interval $I=[a, b]\subset [-\pi, \pi],~a\neq b,$
let $\xi_I$ denote a function with the following properties:
$\xi_I(t)=0$ for $-\pi\leq t\leq a$, $\xi_I(t)=1$ for $b\leq
t\leq \pi$, $\xi_I$ is linear and continuous on $I$. We put
$$
g_n=2^{-n}\sum_{\alpha\in\{0; 1\}^n} \xi_{I_\alpha}.
$$
It is easily seen that the sequence $\{g_n\}$ converges
uniformly on $[-\pi, \pi]$. Let $g=\lim g_n$. We have
$g(-\pi)=0, ~g(\pi)=1$. Let $f=\sin \pi g$. Then $f$ is
constant on each interval complementary to $E$, and
$f\not\equiv\mathrm{const}$.

  Let us show that $f\in H^\omega( T)$. It suffices to
verify that
$$
|g(t_0)-g(t_1)|\leq \mathrm{const}\,\omega(|t_0-t_1|)
\quad\forall t_0, t_1\in E.
$$
Let $t_0, ~t_1\in E$. Then for a certain $k$ we have $t_0,
~t_1\in I_{i_1\ldots i_k}$ and at the same time $t_0\in
I_{i_1\ldots i_k 0}, ~t_1\in I_{i_1\ldots i_k 1}$ (the case
when $t_0\in I_0, ~t_1\in I_1$ can be omitted). Therefore,
$$
|g_n(t_0)-g_n(t_1)|\leq 2^{-k} \quad \forall n\geq k,
$$
whence
$$
|g(t_0)-g(t_1)|\leq 2^{-k}.
$$
By (16) we have $|t_0-t_1|\geq \delta_{k+1}$, so, taking (14)
into account, we obtain
$$
|g(t_0)-g(t_1)|\leq 2\cdot 2^{-(k+1)}\leq
(2/\gamma)\omega(\delta_{k+1})\leq (2/\gamma)\omega(|t_0-t_1|).
$$
The theorem is proved.

\quad

  We note that the function $f$ constructed in the proof
is of bounded variation on $ T$.

  It is not clear for what precisely $\omega$'s the conclusion
of Theorem 2 holds. In the proof we used
$\omega(\delta)=(\log(1/\delta))^{-\alpha}, ~0<\alpha<1$.
Similarly one can show that if
$$
\lim \omega(\delta)\log(1/\delta)=\infty,
$$
then it is impossible to attain (11). Possibly the same holds
for all $\omega$ unless $\omega(\delta)=O(\delta)$.

\begin{center}
\textbf{\S~2. The classes $A_\varepsilon$}
\end{center}

  As we mentioned in Introduction, for every function in
$C(T)$ there is a change of variable which brings it into any
given class $A_\varepsilon(T)$. This result does not extend to
compact families of functions. Moreover there exists a sequence
$\varepsilon_n\rightarrow 0$ such that, in general, there is no
single change of variable which will bring two functions into
$A_\varepsilon$. We do not know what conditions imposed on
$\varepsilon$ are necessary and sufficient in order that for
every pair of (real-valued) continuous functions there is a
change of variable that brings the pair into $A_\varepsilon$.

  From Theorem 1, we immediately obtain the following
corollary.

\quad

\textsc{Corollary 2.} \emph{Let
$$
\sum_{n=1}^\infty \varepsilon_n/n<\infty.
\eqno(17)
$$
Then for every class $H^\omega(T)$ there is a change of
variable which brings it into $A_\varepsilon(T)$, i.e., there
exists a homeomorphism $h$ of the circle $T$ such that $f\circ
h\in A_\varepsilon(T)$ for all $f\in H^\omega(T)$}

\quad

  For the proof it suffices to choose a sequence
$\{\lambda_n\}, ~\lambda_n\rightarrow\infty,$ with
$$
\sum_{n=1}^\infty \lambda_n\varepsilon_n/n<\infty,
$$
and apply Theorem 1.

  We say that a sequence $\varepsilon$ is regular if it
is non-increasing, and $\{n\varepsilon_n\}$ is non-decreasing.

  The following theorem shows that, for regular sequences,
condition (17) is also necessary in order that for every class
$H^\omega(T)$ there is a change of variable which brings
$H^\omega(T)$ into $A_\varepsilon$.

\quad

\textsc{Theorem 3.} \emph{Suppose that the sequence
$\varepsilon$ is regular and
$$
\sum_{n=1}^\infty \varepsilon_n/n=\infty.
\eqno(18)
$$
Then there exists a class $H^\omega(T)$ such that there is no
change of variable which will bring it into $A_\varepsilon(T)$,
i.e., for every homeomorphism $h : T\rightarrow T$ there is an
$f\in H^\omega(T)$ such that $f\circ h\notin
A_\varepsilon(T)$.}

\quad

\textsc{Proof.} We say that a function $f$ is of class
$H_{\mathrm{loc}(0)}^\omega$ if there is an interval
$I\subseteq T$, containing $\{0\}$, such that
$$
\sup_{\underset{t_1, t_2\in I}{|t_1-t_2|<\delta}}
|f(t_1)-f(t_2)|=O(\omega(\delta)), \qquad \delta\rightarrow 0.
$$
It is known that if a function $f\in A( T)$ is monotonic in a
neighborhood of a point, then in a certain neighborhood of this
point the modulus of continuity of $f$ is logarithmic at worst
(Katznelson; see [5], Chapter II, \S~12). Similar result holds
for classes $A_\varepsilon$ provided that $\varepsilon_n$ tends
to zero sufficiently slowly (see [1], Lemma 4.3). The following
lemma shows that this result is valid under the assumption that
(18) is satisfied and $\varepsilon$ is regular.

\quad

\textsc{Lemma 3.} \emph{Under the assumptions of the theorem on
the sequence $\varepsilon$ there exists a function
$\omega_\varepsilon$ satisfying
$\omega_\varepsilon(\delta)\downarrow 0$ as $\delta\downarrow
0$, such that if $f\in A_\varepsilon(T)\cap C(T)$ and $f$ is
monotonic in a neighborhood of zero, then $f\in
H_{\mathrm{loc}(0)}^{\omega_\varepsilon}$.}

\quad

\textsc{Proof.} We identify $T$ with the interval $[-\pi,
\pi]$. Let $\gamma(\theta)$ be the function on $(0, \pi]$ that
takes value $n\varepsilon_n$ for $\theta\in (\pi/(n+1), \pi/n],
~n=1, 2, \ldots$. The function $\gamma(\theta)$ defined in this
way increases as $\theta$ decreases to zero, so
$$
\sup_{0\leq\alpha\leq\pi}\bigg|\int_0^\alpha
\gamma(\theta)\sin k\theta d\theta\bigg|\leq\bigg|\int_0^{\pi/|k|}
\gamma(\theta)\sin k\theta d\theta\bigg|\leq
$$
$$
\leq|k|\bigg|\int_0^{\pi/|k|}
\gamma(\theta)\theta d\theta\bigg|=
|k|\sum_{n\geq |k|}n\varepsilon_n\int_{\pi/(n+1)}^{\pi/n}\theta d\theta
\leq \pi^2\varepsilon_{|k|}, \qquad k\neq 0.
\eqno(19)
$$

\quad

  We now follow the method used to prove Lemma 4.3 in [1].
Note that if $f\in A_\varepsilon( T)\cap C( T)$ and
$0<\delta<\delta_0<\pi$, then
$$
\bigg|\int_{\delta}^{\delta_0}
(f(t+\theta)-f(t-\theta))\gamma(\theta)d\theta\bigg|
\leq 4\pi^2\|f\|_{A_\varepsilon}, \qquad t\in T.
\eqno(20)
$$
Indeed, if $f$ is a trigonometric polynomial, then, using (19),
we obtain
$$
\bigg|\int_{\delta}^{\delta_0}
(f(t+\theta)-f(t-\theta))\gamma(\theta)d\theta\bigg|=
$$
$$
=2\bigg|\sum_{k\in Z}\widehat{f}(k)e^{ikt}
\int_{\delta}^{\delta_0}\gamma(\theta)\sin k\theta d\theta\bigg|
\leq 4\pi^2\|f\|_{A_\varepsilon}.
$$
In the general case one should approximate $f$ by its Fej\'er
sums.

  Let $f$ be monotonic in a $2\delta_0$ -neighborhood of zero.
Then, for every $t$ in the $\delta_0$ -neighborhood of zero and
for every $\delta$ with $0<\delta<\delta_0$ we obtain from (20)
$$
|f(t+\delta)-f(t-\delta)|\int_\delta^{\delta_0}\gamma(\theta) d\theta
\leq\bigg|\int_{\delta}^{\delta_0}
(f(t+\theta)-f(t-\theta))\gamma(\theta)d\theta\bigg|
\leq 4\pi^2\|f\|_{A_\varepsilon}.
$$
Therefore
$$
|f(t+\delta)-f(t)|\int_\delta^{\delta_0}\gamma(\theta) d\theta=O(1)
$$
uniformly with respect to $t, ~|t|<\delta_0$. It remains only
to notice that
$$
\int_\delta^\pi\gamma(\theta) d\theta \uparrow\infty
$$
as $\delta\downarrow 0$, and to put
$$
\omega_\varepsilon(\delta)=
1\bigg/\int_\delta^\pi\gamma(\theta) d\theta.
$$
The lemma is proved.

\quad

  Let us now show that the smoothness of all functions
in $H^\omega(T)$ which are monotonic in a neighborhood of zero
cannot be improved even locally.

\quad

\textsc{Lemma 4.} \emph{Let
$\omega(2\delta)\leq2\omega(\delta)$ for all $\delta>0$. Let
$\overline{\omega}(\delta)=o(\omega(\delta))$ as
$\delta\rightarrow 0$, and let $h$ be a homeomorphism $
T\rightarrow T$ with $h(0)=0$. Then there exists a function
$f\in H^\omega( T)$, monotonic in a neighborhood of zero, such
that $f\circ h\notin H_{\mathrm{loc}(0)}^{\overline{\omega}}$.}

\quad

\textsc{Proof.} Similarly to what we did to prove Theorem 2,
for each closed interval $[a, b]\subset [-\pi, \pi]$ consider a
function $\xi_I$ continuous on $[-\pi, \pi]$, such that
$\xi_I(t)=0$ for $-\pi\leq t\leq a$, $\xi_I(t)=1$ for $b\leq
t\leq \pi$, and $\xi_I$ is linear on $I$.

  We fix a positive sequence
$\lambda_n, ~n=1, 2, \ldots, ~\lambda_n\rightarrow 0,$ such
that
$$
\frac{\overline{\omega}(1/n)}{\omega(\lambda_n/n)}\rightarrow 0,
\qquad n\rightarrow\infty.
$$
Choose a sequence of positive integers $\{n_k\}$ with the
following properties
$$
\sum_{k=1}^\infty\lambda_{n_k}<\pi/20,
$$
$$
\omega(\lambda_{n_{k+1}}/n_{k+1})\leq\frac{1}{2}\omega(\lambda_{n_k}/n_k),
\qquad k=1, 2, \ldots.
$$
Let $d_k=\lambda_{n_k}/n_k$. Then
$$
\overline{\omega}(1/n_k)=o(\omega(d_k)), \qquad k\rightarrow\infty,
\eqno(21)
$$
$$
\sum_{k=1}^\infty 20n_k d_k=d<\pi,
\eqno(22)
$$
$$
\omega(d_{k+1})\leq\frac{1}{2}\omega(d_k), \qquad k=1, 2, \ldots.
\eqno(23)
$$

   Choose points $a_k, ~k=1, 2,
\ldots,$ in the interval $[0, d]$ so that $0<a_{k+1}<a_k,
~a_k\rightarrow 0$, and $|a_k-a_{k+1}|=20n_kd_k$; this is
possible since (22). We partition each interval $[a_{k+1},
a_k]$ into $20n_k$ closed intervals of length $d_k$.

  Let $h$ be a homeomorphism $ T\rightarrow T$
with $h(0)=0$. For $k=1, 2, \ldots$ we exclude from the
intervals that form the partition of $[a_{k+1}, a_k]$ the
left-hand one and find among the remaining ones an interval
$I_k$ such that
$$
|h^{-1}(I_k)|\leq 2\pi/(20n_k-1).
\eqno(24)
$$
By our construction, if $k_1<k_2$, then there is an interval of
length $d_{k_1}$ between $I_{k_1}$ and $I_{k_2}$.

  We put
$$
g=\sum_{k=1}^\infty \omega(d_k)\xi_{I_k}
$$
and define a function $f\in C( T)$ as follows: $f(t)=g(t)$ for
$-\pi\leq t\leq d, ~f(\pi)=0$, and $f$ is linear on $[d, \pi]$.
Obviously $f$ is monotonic in a neighborhood of zero.

  Let us show that $f\in H^\omega( T)$. It suffices to
verify that
$$
|g(t_1)-g(t_2)|\leq\mathrm{const}\cdot\omega(|t_1-t_2|)
$$
for $t_1, t_2\in\bigcup_k I_k$. Assume first that $t_1, t_2\in
I_k, ~t_1\neq t_2$. Then, choosing a positive integer $p$ such
that
$$
2^p\leq\frac{d_k}{|t_1-t_2|}<2^{p+1},
$$
we have
$$
|g(t_1)-g(t_2)|=\frac{\omega(d_k)}{d_k}|t_1-t_2|
\leq\omega(2^{p+1}|t_1-t_2|)2^{-p}
\leq
$$
$$
\leq 2^{p+1}\omega(|t_1-t_2|)2^{-p}=2\omega(|t_1-t_2|).
$$
Assume now that $t_1\in I_{k_1}$ and $t_2\in I_{k_2},
~k_1<k_2$. Then $|t_1-t_2|\geq d_{k_1}$, and it follows from
(23) that
$$
|g(t_1)-g(t_2)|\leq\sum_{k_1\leq k\leq k_2}\omega(d_k)
\leq 2\omega(d_{k_1})\leq2\omega(|t_1-t_2|).
$$

  Let us show that
$f\circ h\notin H_{\mathrm{loc}(0)}^{\overline{\omega}}$.
Assuming the contrary, from (24) we obtain
$$
\omega(d_k)=|f(I_k)|=|f\circ h(h^{-1}(I_k))|=
O(\overline{\omega}(|h^{-1}(I_k)|))=
O(\overline{\omega}(1/n_k)), \quad k\rightarrow\infty,
$$
which contradicts (21). The lemma is proved.

\quad

  We shall now complete the proof of the theorem. Let $\omega$
be an increasing continuous function on $[0, \infty)$ with
$\omega(0)=0, ~\omega(2\delta)\leq 2\omega(\delta)
~\forall\delta>0,$ and
$\omega_\varepsilon(\delta)=o(\omega(\delta))$, where
$\omega_\varepsilon$ is the function from Lemma 3. For example
these conditions hold for
$$
\omega(\delta)=\sup_{|t_1-t_2|\leq
\delta}|\sqrt{\omega_\varepsilon(t_1)}-\sqrt{\omega_\varepsilon(t_2)}|.
$$
Suppose that for some homeomorphism $h$ of $ T$ onto itself we
have $f\circ h\in A_\varepsilon(T)$ for all $f\in H^\omega(
T)$. We may assume that $h(0)=0$. By Lemma 4, there exists a
function $f\in H^\omega( T)$, monotonic in a neighborhood of
zero, such that $f\circ h\notin
H_{\mathrm{loc}(0)}^{\omega_\varepsilon}$. But, since $f\circ
h$ is monotonic in a neighborhood of zero, it follows from
Lemma 3 that $f\circ h\in
H_{\mathrm{loc}(0)}^{\omega_\varepsilon}$. The contradiction
proves the theorem.

\begin{center}
\textbf{\S~3. Sobolev classes}
\end{center}

  Let $W_2^\lambda(T)$ be the class of all functions $x$ with
$$
\|x\|_{W_2^\lambda}=\bigg(\sum_{k\in Z}
(|\widehat{x}(k)||k|^\lambda)^2\bigg)^{1/2}<\infty.
$$

\quad

 Theorem 1 implies the following corollary:

\quad

\textsc{Corollary 3.} \emph{For every class $H^\omega(T)$ there
exists a change of variable which brings it into
$\bigcap_{\lambda<1/2}W_2^\lambda(T)$, i.e., there exists a
homeomorphism $h$ of the circle $T$ onto itself such that
$f\circ h\in\bigcap_{\lambda<1/2}W_2^\lambda(T)$ for all $f\in
H^\omega(T)$.}

\quad

  Let us recall that for every function in $C(T)$ there is a change of
variable which brings it into $W_2^{1/2}$. Is it true that for
every class $H^\omega$ there exists a change of variable which
brings it into $W_2^{1/2}$? The answer to this question posed
in [3, p.41] is negative. Moreover the following theorem holds.

\quad

\textsc{Theorem 4.} \emph{Let $f\in C( T)$. The following
conditions are equivalent:}

(i) \emph{For every function $g\in C(T)$ there exists a change
of variable which brings the pair $\{f, g\}$ into
$W_2^{1/2}(T)$, i.e., there is a homeomorphism $h :
T\rightarrow T$ such that $f\circ h\in W_2^{1/2}$ and $g\circ
h\in W_2^{1/2}$ .}

(ii) \emph{f is of bounded variation on $ T$.}

\quad

\textsc{Proof.} Note that the obvious estimate
$$
c_1|k|\leq \int_0^1\bigg(\frac{\sin k\delta}{\delta}\bigg)^2 d\delta
\leq c_2|k|, \qquad k\in Z
$$
(where $c_1, c_2>0$ are independent of $k$) and the identity
$$
\frac{1}{2\pi}\int_ T
|x(t+\delta)-x(t-\delta)|^2 dt=
4\sum_{k\in Z}|\widehat{x}(k)|^2\sin^2k\delta
$$
imply the equivalence of the seminorms $\|\cdot\|_{W_2^{1/2}}$
and $\|\cdot\|$, where
$$
\|x\|=\bigg(\int_0^1\frac{1}{\delta^2}\int_ T
|x(t+\delta)-x(t-\delta)|^2 dt d\delta\bigg)^{1/2}.
$$
We shall use this fact later.

By $\mathrm{Var}(x, E)$ we denote the variation of a function
$x(t)$ on a set $E\subseteq T$.

  1) $\mathrm{(i)\Rightarrow (ii)}$. Let $f\in C( T)$.
Under the assumption that $\mathrm{Var}(f, T)=\infty$ we shall
construct a function $g\in C( T)$ such that there is no single
change of variable which will bring both $f$ and $g$ into
$W_2^{1/2}$. Thus the implication $\mathrm{(i)\Rightarrow
(ii)}$ will be prowed.

\quad

\textsc{Lemma 5.} \emph{There exists a monotonic sequence
$t_k\in T, ~k=1, 2, \ldots,$ such that
$$
\sum_{k=1}^\infty |f(t_{k+1})-f(t_k)|=\infty.
$$}

\quad

\textsc{Proof.} Note that there exists a point $\theta\in T$
such that $f$ has infinite variation in every neighborhood of
$\theta$. Indeed, otherwise each point $t\in T$ would have a
neighbourhood $U_t$ such that $\mathrm{Var}(f, U_t)<\infty$,
and choosing a finite covering of the circle from the family
$\{U_t, ~t\in T\}$ we would obtain $\mathrm{Var}(f, T)<\infty$
which contradicts the assumption.

  Fix $\theta\in T$ with the indicated property. Then
either $\mathrm{Var}(f, (\theta', \theta))=\infty$ for every
open interval $(\theta', \theta), ~\theta'<\theta,$ or
$\mathrm{Var}(f, (\theta, \theta'))=\infty$ for every open
interval $(\theta, \theta'), ~\theta'>\theta$.

  Consider the first case (the second one is similar).
Since $\mathrm{Var}(f, (-\pi, \theta))=\infty$, one can find
points $t_k, ~k=1, \ldots, n_1,
~-\pi<t_1<\ldots<t_k<t_{k+1}<\ldots<t_{n_1}<\theta$, such that
$$
\sum_{k=1}^{n_1-1}|f(t_{k+1})-f(t_k)|>1.
$$
Assume that the points $t_1<\ldots t_{n_s}<\theta$ have already
been defined. Since $\mathrm{Var}(f, (t_{n_s},
\theta))=\infty$, one can find points $t_k, ~k=n_s+1, \ldots,
n_{s+1}, ~t_{n_s}<t_{n_s+1}<\ldots<t_{n_{s+1}}<\theta$ such
that
$$
\sum_{k=n_s+1}^{n_{s+1}-1}|f(t_{k+1})-f(t_k)|>1.
$$
Continuing this process, we obtain the required sequence.

\quad

\textsc{Lemma 6.} \emph{There exist a function $g\in C( T)$ and
a sequence of functions $g_n, ~n=1, 2, \ldots,$ with the
following properties:
$$
|g_n(t_1)-g_n(t_2)|\leq |g(t_1)-g(t_2)| \qquad \forall t_1, t_2,
\eqno(25)
$$
$$
\mathrm{Var}(g_n,  T)<\infty, \quad \forall n,
\eqno(26)
$$
$$
\sup_n\bigg|\int_ T f(t)dg_n(t)\bigg|=\infty.
\eqno(27)
$$
}

\quad

\textsc{Proof.} Let $\{t_k\}$ be the sequence from Lemma 5. Put
$M^+=\{k : f(t_{k+1})-f(t_k)>0\}, ~M^-=\{k :
f(t_{k+1})-f(t_k)<0\}$. Then at least one of the sums
$$
\sum_{k\in M^+}|f(t_{k+1})-f(t_k)|, \qquad
\sum_{k\in M^-}|f(t_{k+1})-f(t_k)|
$$
is infinite. We proceed with our construction under the
assumptions that $t_1<t_2<\ldots<t_k<t_{k+1}<\ldots$ and
$$
\sum_{k\in M^+}|f(t_{k+1})-f(t_k)|=\infty
$$
(the other cases are similar).

  Let $k_j, ~j=1, 2, \ldots,$ be the sequence of numbers that
form $M^+$. Put $a_j=t_{k_j}, ~b_j=t_{k_j+1}, ~j=1, 2, \ldots$.
The intervals $(a_j, b_j)$ are pairwise disjoint and
$$
f(b_j)-f(a_j)>0, \qquad \sum_j(f(b_j)-f(a_j))=\infty.
\eqno(28)
$$

\quad

  Choose a sequence of positive numbers $\gamma_j, ~j=1, 2,
\ldots,$ that decreases to zero and satisfies
$$
\sum_j\gamma_j(f(b_j)-f(a_j))=\infty.
\eqno(29)
$$
Choose also a sequence $\{\varepsilon_j\}$ such that
$0<2\varepsilon_j<b_j-a_j$ and
$$
\omega(\varepsilon_j, f)<2^{-j}.
\eqno(30)
$$

  For $j=1, 2, \ldots$ we define functions $\lambda_j$ as
follows: $\lambda_j\in C( T)$, $\lambda_j(t)=0$ for $t\notin
(a_j, b_j)$, $\lambda_j(t)=\gamma_j$ for $t\in
(a_j+\varepsilon_j, b_j-\varepsilon_j)$, and $\lambda_j$ is
linear on $(a_j, a_j+\varepsilon_j)$ and $(b_j-\varepsilon_j,
b_j)$.

  Let
$$
g(t)=\sum_{j=1}^\infty\lambda_j(t)
$$

  Obviously $g\in C( T)$. Let
$$
g_n(t)=\max\{\gamma_n, g(t)\}.
$$
It is easily verified that (25) and (26) hold.

  Let us show that (27) holds. Note that the set
$\{t : \gamma_n<\lambda_j(t)<\gamma_j\}$ is empty if $j\geq n$,
whereas if $1\leq j<n$ it consists of two intervals $I_{j n}^+$
and $I_{j n}^-$ of equal length on which $g$ increases and
decreases respectively. Note also that
$$
|I_{j n}^{\pm}|\leq\varepsilon_j,
\eqno(31)
$$
$$
\lim_n |I_{j n}^{\pm}|=\varepsilon_j \qquad \forall j.
\eqno(32)
$$

  We have
$$
\int_ T f(t)dg_n(t)=\sum_{j=1}^{n-1}
\bigg(\int_{I_{j n}^+}f(t)\frac{\gamma_j}{\varepsilon_j}dt-
\int_{I_{j n}^-}f(t)\frac{\gamma_j}{\varepsilon_j}dt\bigg)=
$$
$$
=\sum_{j=1}^{n-1}\frac{\gamma_j}{\varepsilon_j}
\bigg[\int_{I_{j n}^+}(f(t)-f(a_j))dt+f(a_j)|I_{j n}^+|-
\int_{I_{j n}^-}(f(t)-f(b_j))dt-f(b_j)|I_{j n}^-|\bigg].
$$
Taking (31) and (30) into account, we obtain
$$
\bigg|\int_ T f(t)dg_n(t)+\sum_{j=1}^{n-1}
\gamma_j(f(b_j)-f(a_j))\frac{|I_{j n}^\pm|}{\varepsilon_j}\bigg|\leq
$$
$$
\leq\sum_{j=1}^{n-1}\frac{\gamma_j}{\varepsilon_j}
\bigg(\int_{I_{j n}^+}|f(t)-f(a_j)|dt+
\int_{I_{j n}^-}|f(t)-f(b_j)|dt\bigg)\leq
$$
$$
\leq\sum_{j=1}^{n-1}2\gamma_j\omega(\varepsilon_j, f)=O(1),
\qquad n\rightarrow\infty.
\eqno(33)
$$
Suppose that (27) does not hold. Then (see (33))
$$
\sum_{j=1}^{n-1}
\gamma_j(f(b_j)-f(a_j))\frac{|I_{j n}^\pm|}{\varepsilon_j}=O(1).
\eqno(34)
$$
Since the terms in (34) are positive (see (28)), we have for
$m<n$
$$
\sum_{j=1}^m
\gamma_j(f(b_j)-f(a_j))\frac{|I_{j n}^\pm|}{\varepsilon_j}
\leq\mathrm{const}.
\eqno(35)
$$
Using (32) and taking the limit in (35), we obtain
$$
\sum_{j=1}^m
\gamma_j(f(b_j)-f(a_j))\leq\mathrm{const}.
$$
Since $m$ is arbitrary, this contradicts (29). The lemma is
proved.

\quad

Let $g$ be a function as in Lemma 6. Let us show that there is
no change of variable which will bring the pair $\{f, g\}$ into
$W_2^{1/2}$. Suppose that, on the contrary, $f\circ h\in
W_2^{1/2}$ and $g\circ h\in W_2^{1/2}$ for a certain
homeomorphism $h : T\rightarrow T$. Then using (25) and the
equivalence of the seminorms $\|\cdot\|_{W_2^{1/2}}$ and
$\|\cdot\|$ we see that $g_n\circ h\in W_2^{1/2}$ for all $n$,
and
$$
\|g_n\circ h\|_{W_2^{1/2}}\leq\mathrm{const}\cdot\|g\circ h\|_{W_2^{1/2}}.
\eqno(36)
$$
Note now, that if $x, y\in W_2^{1/2}( T)\cap C( T)$ and $y$ is
a function of bounded variation, then
$$
\frac{1}{2\pi}\bigg|\int_{ T}x(t)dy(t)\bigg|\leq
\|x\|_{W_2^{1/2}}\|y\|_{W_2^{1/2}}.
$$
(This is obvious if $x$ is a trigonometric polynomial; in the
general case one should approximate $x$ by Fej\'er sums.) Thus,
from (36), taking (26) into account, we have
$$
\frac{1}{2\pi}\bigg|\int_{ T}f(t)dg_n(t)\bigg|=
\frac{1}{2\pi}\bigg|\int_{ T}f\circ h(t)dg_n\circ h(t)\bigg|
\leq \|f\circ h\|_{W_2^{1/2}}\|g_n\circ h\|_{W_2^{1/2}}=O(1).
$$
This contradicts (27). The implication (i)$\Rightarrow$(ii) is
proved.

2) (ii)$\Rightarrow$(i). The method of the proof of Theorem A
that uses a conformal mapping of the disk onto a suitable
domain (see [1], Proof of Theorem 3.1) admits the following
modification that allows to bring an arbitrarily given function
$g\in C( T)$ into $W_2^{1/2}$. Assuming that $g(t)\geq
\gamma>0$ for all $t$ (this does not restrict generality) we
put $Q(t)=g(t)e^{it}$. Obviously, when $t$ runs over $[-\pi,
\pi]$ the point  $Q(t)$ describes a simple closed curve
$\Gamma$ in the complex plane. Let $D$ be a domain bounded by
this curve, and let $\Phi$ be a conformal mapping of the disk
$|z|<1$ onto $D$. Then $\Phi$ extends continuously to the
circle $|z|=1$, and provides a homeomorphism of the circle onto
$\Gamma$. Thus the function $\varphi(t)=\Phi(e^{it})$ has the
form $\varphi=Q\circ h$ where $h$ is a homeomorphism $
T\rightarrow T$.

  It is well known that $\pi\sum_{n\geq 0} |\widehat{\varphi}(n)|^2
n$ is the area of $D$. Therefore $Q\circ h\in W_2^{1/2}(T)$.
Using the equivalence of the seminorms $\|\cdot\|_{W_2^{1/2}}$
and $\|\cdot\|$, we obtain $|Q\circ h|\in W_2^{1/2}$ and it
remains to note that $|Q\circ h|=g\circ h$. This modification
of the P\'al--Bohr theorem was found by A. M. Olevski\v{\i}
(personal communication).

   We now note that if $x\in W_2^{1/2}\cap C( T)$ and
$x(t)\geq\gamma>0$ for all $t$, then $1/x\in W_2^{1/2}$; this
follows from the equivalence of the seminorms. This equivalence
implies also that if two continuous complex-valued functions
are of class $W_2^{1/2}$ then their product is in $W_2^{1/2}$.
So, in the construction described above, we obtain in addition
that
$$
e^{ih}=\frac{1}{|Q\circ h|}Q\circ h\in W_2^{1/2}.
$$

  We have thus proved the following lemma.

\quad

\textsc{Lemma 7.} \emph{For every pair $\{g(t), e^{it}\}$,
where $g\in C(T)$, there is a change of variable which brings
it into $W_2^{1/2}(T)$, i.e., there exists a homeomorphism $h :
T\rightarrow T$ such that $g\circ h\in W_2^{1/2}(T)$ and
$e^{ih}\in W_2^{1/2}(T)$.}

\quad

  We now complete the proof of the implication
(ii)$\Rightarrow$(i). Suppose that (ii) holds and $g$ is an
arbitrary function in $C( T)$. It is easily seen that a
function, which is continuous and of bounded variation on the
circle, can be turned into a Lipschitz function by an
appropriate change of variable. Fix a homeomorphism $\psi :
T\rightarrow T$ such that
$$
|f\circ\psi(t_1)-f\circ\psi(t_2)|\leq
\mathrm{const}\cdot |e^{it_1}-e^{it_2}| \qquad \forall t_1, ~t_2.
$$

  Applying Lemma 7 to the function $g\circ\psi$, we obtain a
homeomorphism $h : T\rightarrow T$ such that $g\circ\psi\circ
h\in W_2^{1/2}(T)$ and
$$
e^{ih}\in W_2^{1/2}(T).
\eqno(37)
$$

  It remains only to observe that, taking account of the
equivalence of the seminorms $\|\cdot\|_{W_2^{1/2}}$ and
$\|\cdot\|$, from (37) and the inequality
$$
|f\circ\psi\circ h(t_1)-f\circ\psi\circ h(t_2)|\leq
\mathrm{const}\cdot |e^{ih(t_1)}-e^{ih(t_2)}| \qquad \forall t_1, ~t_2
$$
it follows that $f\circ\psi\circ h\in W_2^{1/2}( T)$.

\quad

  The author thanks A. M. Olevski\v{\i} for his help and
attention.

\quad

\begin{center}
\textbf{References}
\end{center}

\flushleft
\begin{enumerate}

\item Olevski\v{\i}, A. M., Modifications of functions and
    Fourier series, \emph{Uspekhi Mat. Nauk}, \textbf{40}:3
    (1985), 157--193 (Russian); English transl.:
    \emph{Russian Math. Surveys}, \textbf{40}:3 (1985),
    181--224.

\item Olevski\v{\i}, A. M., Homeomorphisms of the circle,
    modifications of functions, and Fourier series,
    \emph{Proc. Int. Congr. Math.} (Berkeley Calif., 1986),
    Vol. 2, Amer. Math. Soc., Providence, R. I., 1987, pp.
    976--989; English transl.: \emph{Amer. Math. Soc.
    Transl.} (2) \textbf{147} (1990).

\item Olevski\v{\i}, A. M., Modifications of functions and
    Fourier series, \emph{Theory of Functions and
    Approximations, Proc. 2nd Saratov Winter School}
    (Saratov, 1984), Part 1, Saratov University, Saratov,
    1986, pp. 31--43 (Russian).

\item Bari, N. K., \emph{Trigonometric Series}, Fizmatgiz,
    Moscow, 1961 (Russian); English transl.: Vols I, II,
    Pergamon Press, Oxford, and Macmillan, New York, 1964.

\item Kahane, J.-P., \emph{S\'eries de Fourier absolument
    convergantes}, Springer-Verlag, Berlin--Heidelberg--New
    York, 1970.

\end{enumerate}

\quad

\quad

Moscow Institute of Electronics and Mathematics,

National Research University Higher School of Economics

E-mail address: \emph{lebedevhome@gmail.com}

\end{document}